\documentclass{siamltex}

\usepackage[utf8]{inputenc}
\usepackage{geometry}

\usepackage{amsmath, amssymb, amsfonts}

\usepackage{graphicx}
\usepackage{xcolor}
\usepackage{tikz}
\usepackage{bm}

\usepackage{algorithm}
\usepackage{algorithmicx}
\usepackage{algpseudocode}

\usepackage{booktabs, multirow, colortbl}

\usepackage{subfig}

\usepackage{enumitem}
\usepackage{titlesec}

\usepackage[colorlinks=true,linkcolor=blue,citecolor=blue,urlcolor=blue]{hyperref}

\let\cite\relax
\usepackage[numbers,sort&compress]{natbib}

\providecommand{\qed}{\hfill\ensuremath{\square}}  

\newtheorem{remark}[theorem]{Remark}
\newtheorem{assumption}[theorem]{Assumption}

\title{Efficient Krylov--Regularization Solvers for Multiquadric RBF Discretizations of the 3D Helmholtz Equation}
\author{
  M. El Guide\thanks{FGSES, University Mohammed VI Polytechnic, Rabat, Morocco}
  \and K. Jbilou\thanks{LMPA, University Littoral Côte d’Opale, Calais, France}
  \and K. Lachhab\thanks{FST, University Cadi Ayyad, Marrakech, Morocco}
  \and D. Ouazar\thanks{UM5R, EMI, Rabat, Morocco}
}
\date{}

\begin{document}
 \maketitle
\begin{abstract}
Meshless collocation with multiquadric radial basis functions (MQ-RBFs) delivers high accuracy for the 3D Helmholtz equation but produces dense, severely ill-conditioned linear systems. We develop and evaluate three complementary methods that embed regularization in Krylov projections to overcome this instability at scale: (i) an inexpensive TSVD  that replaces the full SVD by a short Golub–Kahan bidiagonalization and a tiny projected SVD, retaining the dominant spectral content at greatly reduced cost; (ii) Classical Tikhonov regularization with principled parameter choice (GCV/L-curve), expressed in SVD form for transparent filtering; and (iii) a Hybrid Krylov-Tikhonov (HKT) scheme that first projects with Golub-Kahan and then selects the regularization parameter on the reduced problem, yielding stable solutions in few iterations. Extensive tests on canonical domains (cube and sphere) and a realistic industrial pump-casing geometry demonstrate that HKT consistently matches or surpasses the accuracy of full TSVD/Tikhonov at a fraction of the runtime and memory, while  inexpensive TSVD provides the fastest viable reconstructions when only the leading modes are needed. These results show that coupling Krylov projection with TSVD/Tikhonov regularization provides a robust, scalable pathway for MQ-RBFs Helmholtz methods in complex 3D settings.
\end{abstract}

\begin{keywords}
Helmholtz equation, Golub--Kahan bidiagonalization, Meshless methods, Radial basis function, Tikhonov regularization, Krylov subspaces
\end{keywords}

\section{Introduction}

The accurate and efficient numerical solution of partial differential equations (PDEs) remains a central challenge in computational science and engineering. Among these, the Helmholtz equation occupies a fundamental position due to its role in modeling steady-state wave phenomena in acoustics, electromagnetics, and elasticity. Over the past decades, classical discretization techniques—finite element, finite difference, and finite volume methods—have established themselves as reliable and mathematically rigorous tools for solving a broad range of PDEs \cite{strang1973analysis,brenner2008mathematical,leveque2002finite,eymard2000finite}. However, their reliance on structured or unstructured meshes poses a significant obstacle in the context of complex geometries, multi-scale domains, and three-dimensional problems, where mesh generation and adaptation can dominate both the modeling effort and the computational cost.

Meshless methods based on radial basis functions (RBFs) have emerged as an attractive alternative that circumvents the need for explicit mesh connectivity. By constructing the approximation directly from scattered nodes, these methods offer remarkable geometric flexibility, smoothness, and the potential for spectral-like accuracy \cite{fornberg2015primer,Sarra2009}. Among the wide family of RBFs, the multiquadric (MQ) function introduced by Hardy and later extended by Kansa has demonstrated superior approximation capabilities and robustness for elliptic and wave equations \cite{Kansa1990,Sarra2009}. The Kansa collocation approach, in particular, provides a conceptually simple and general framework for discretizing PDEs on arbitrary geometries, including curved or multiply connected domains.

Despite these appealing properties, global RBF formulations suffer from two major numerical difficulties that limit their scalability and practical applicability. First, the use of globally supported basis functions results in dense algebraic systems whose computational and memory requirements grow quadratically with the number of nodes. Second, and more critically, the collocation matrices become increasingly ill-conditioned as the RBF shape parameter $\varepsilon$ decreases, a phenomenon that reflects a fundamental trade-off between approximation accuracy and numerical stability. Small values of $\varepsilon$ improve the approximation power but lead to nearly linearly dependent columns in the interpolation matrix, causing the condition number to grow rapidly and rendering the system extremely sensitive to round-off errors \cite{schaback1995error,fornberg2011stable,fasshauer2015kernel}. This intrinsic instability—often referred to as the RBF uncertainty principle—has motivated the development of stabilization techniques such as regularization, spectral filtering, and localized approximations \cite{Sarra2009,Hansen1998}.

The present work proposes an efficient and stable computational strategy that directly addresses the ill-conditioning of MQ-based discretizations while retaining their superior approximation properties. Three complementary algorithms are developed within a unified framework.  
First, an \emph{Inexpensive Truncated Singular Value Decomposition (Ine-TSVD)} is introduced, which replaces the full SVD with a low-cost projection based on the Golub–Kahan bidiagonalization process. This projection isolates the dominant spectral components of the RBF system at a fraction of the computational cost, enabling stable approximations even for large three-dimensional problems.  
Second, a classical \emph{Tikhonov regularization} scheme is incorporated, with the regularization parameter selected by data-driven strategies such as the generalized cross-validation (GCV) and L-curve criteria \cite{Tikhonov1977,Golub1979,Hansen1993,Hansen1998}.  
Finally, a \emph{Hybrid
Krylov–Tikhonov method} combines both ideas: the regularization parameter is estimated on a reduced Krylov subspace, where the dominant subspace is captured through partial bidiagonalization. This hybrid approach achieves a favorable balance between computational efficiency, stability, and accuracy.

The proposed methods transform the ill-conditioned dense systems arising from MQ collocation into well-posed reduced problems, enabling accurate and scalable computation for three-dimensional Helmholtz equations. By exploiting the spectral filtering effect of truncated and regularized projections, the methods significantly improve robustness without compromising accuracy. Numerical experiments on canonical test domains (unit cube and unit sphere) and a realistic industrial pump-casing geometry demonstrate that the Hybrid
Krylov–Tikhonov method provides an excellent compromise between accuracy and computational cost, while the inexpensive truncated SVD offers a particularly efficient option when only the dominant modes of the solution are sought.

The remainder of the paper is organized as follows. Section~\ref{sec:MQ} reviews the multiquadric RBF framework for the three-dimensional Helmholtz equation. Section~\ref{sec:krylov} introduces the proposed regularization techniques combined with Krylov subspace methods. Section~\ref{sec:existence-convergence} presents theoretical results on existence, uniqueness, and convergence. Numerical experiments are reported in Section~\ref{sec:experiments}, and concluding remarks are given in Section~\ref{sec:conclusion}.

\section{The Multiquadric Radial Basis Function Method (MQ-RBF)}
\label{sec:MQ}

\subsection{The MQ-RBF Collocation Method}
We consider the three-dimensional acoustic Helmholtz equation defined on a bounded domain $\Omega \subset \mathbb{R}^3$:
\begin{equation}
\Delta u(X) + k^2\,u(X) = f(X), \qquad X \in \Omega, 
\label{eq:Helmholtz}
\end{equation}
supplemented with boundary conditions of Robin type:
\begin{equation}
a(X)\,u(X) + b(X)\,\frac{\partial u}{\partial n}(X) = g(X), \qquad X \in \partial\Omega,
\label{eq:Robin}
\end{equation}
where $\partial u/\partial n$ denotes the outward normal derivative on the boundary, and $a(X), b(X), g(X)$ are given functions. In this work, we focus on Dirichlet boundary conditions (taking $a(X)\neq 0$, $b(X)=0$, so that \eqref{eq:Robin} reduces to $u(X) = g(X)$ on $\partial\Omega$).

The MQ-RBF, originally introduced by Hardy in the context of cartography \cite{Hardy1971} and later established as a powerful interpolation tool by Franke \cite{Franke1982} and Micchelli \cite{Micchelli1986}, is defined as 
\begin{equation}
\phi(r;\varepsilon) = \sqrt{\,1 + \varepsilon^2 r^2\,}\,,
\qquad r = \|X - X_j^c\|_2,
\label{eq:MQ}
\end{equation}
where $\varepsilon > 0$ is the shape parameter controlling the flatness of the basis, and $\{X_j^c\}_{j=1}^N \subset \mathbb{R}^3$ are the center points. The MQ is globally supported, infinitely smooth, and non-polynomial, making it particularly effective for high-accuracy approximation of smooth fields \cite{Sarra2009}. Kansa’s pioneering work \cite{Kansa1990} first applied MQ-RBFs to PDEs, sparking widespread interest and further theoretical development \cite{Madych1990}.

\textbf{Collocation framework.} Let the set of collocation centers be partitioned into:
\begin{itemize}
    \item \emph{Interior nodes}: $\{X_i^I\}_{i=1}^{N_I} \subset \Omega$, where the PDE \eqref{eq:Helmholtz} is enforced.
    \item \emph{Boundary nodes}: $\{X_i^B\}_{i=1}^{N_B} \subset \partial\Omega$, where boundary conditions \eqref{eq:Robin} are imposed.
\end{itemize}
The total number of centers is $N = N_I + N_B$. We approximate the solution $u(X)$ as a linear combination of shifted MQ basis functions:
\begin{equation}
u(X) \;\approx\; \sum_{j=1}^N \alpha^{(j)}\, \phi(\|X - X_j^c\|_2;\varepsilon)\,,
\label{eq:u_expansion}
\end{equation}
where $\boldsymbol{\alpha} = [\alpha^{(1)}, \ldots, \alpha^{(N)}]^T \in \mathbb{R}^N$ are unknown coefficients. The interpolation (collocation) matrix $B \in \mathbb{R}^{N\times N}$ associated with the basis functions is defined by
\begin{equation}
B_{ij} = \phi(\,\|X_i^c - X_j^c\|_2;\varepsilon\,)\,, \qquad 1 \le i,j \le N,
\label{eq:B_matrix}
\end{equation}
so that, in matrix form, 
\begin{equation}
\mathbf{u} \approx B\,\boldsymbol{\alpha}\,,
\label{eq:u_matrix}
\end{equation}
where $\mathbf{u}$ is the vector of approximate solution values at the collocation centers.

Applying the differential operator of the Helmholtz equation to the MQ expansion \eqref{eq:u_expansion} yields the interior collocation equations. In particular, for each interior node $X_i^I$, 
\begin{equation}
\mathcal{L}\,\phi(\|X_i^I - X_j^c\|_2;\varepsilon) \;:=\; \Delta \phi(\|X_i^I - X_j^c\|_2;\varepsilon) + k^2\,\phi(\|X_i^I - X_j^c\|_2;\varepsilon)\,,
\label{eq:Lphi}
\end{equation}
for $j = 1,\ldots,N$. This defines the differentiation matrix $H \in \mathbb{R}^{N_I \times N}$ with entries
\begin{equation}
H_{ij} = \Delta \phi(\|X_i^I - X_j^c\|_2;\varepsilon) + k^2\,\phi(\|X_i^I - X_j^c\|_2;\varepsilon)\,,
\label{eq:H_entries}
\end{equation}
for $1 \le i \le N_I$, $1 \le j \le N$. The interior equations can thus be expressed as 
\begin{equation}
H\,\boldsymbol{\alpha} = \mathbf{f}_I\,,
\label{eq:interior_system}
\end{equation}
where $\mathbf{f}_I \in \mathbb{R}^{N_I}$ contains the values of the source term $f(X)$ at interior collocation points. 

For the boundary conditions \eqref{eq:Robin}, let $R \in \mathbb{R}^{N_B \times N}$ denote the discretization (collocation) matrix for the boundary operator, and let $\mathbf{g}_B \in \mathbb{R}^{N_B}$ be the vector of boundary data. The full collocation system combining interior and boundary conditions can be written in block form as 
\begin{equation}
\begin{pmatrix} H \\[1mm] R \end{pmatrix} \boldsymbol{\alpha} \;=\; \begin{pmatrix} \mathbf{f}_I \\[1mm] \mathbf{g}_B \end{pmatrix}, 
\label{eq:augmented_system}
\end{equation}
or equivalently 
\begin{equation}
A\,\boldsymbol{\alpha} = \mathbf{f}, \qquad A \in \mathbb{R}^{N\times N},
\label{eq:full_system}
\end{equation}
where $A$ is the fully populated (dense) collocation matrix and $\mathbf{f}$ is the right-hand side assembled from $\mathbf{f}_I$ and $\mathbf{g}_B$.

\textbf{Conditioning issues.} The conditioning of the system \eqref{eq:full_system} is highly sensitive to the shape parameter $\varepsilon$. As $\varepsilon \to 0$ (i.e., as the basis functions become increasingly flat), the columns of $A$ tend toward linear dependence, causing the condition number $\kappa(A)$ to grow rapidly. This phenomenon reflects the inherent uncertainty principle of RBF methods \cite{Sarra2009}: increasing accuracy through smaller shape parameters inevitably exacerbates ill-conditioning. Consequently, the system becomes numerically unstable and highly sensitive to perturbations.

A deeper understanding of this instability can be obtained by examining the singular value decomposition (SVD) of the differentiation matrix $A$:
\begin{equation}
A = U\,\Sigma\,V^{T},
\label{eq:svdH}
\end{equation}
where $\Sigma = \mathrm{diag}(\sigma_1,\ldots,\sigma_N)$ is the diagonal matrix of singular values satisfying 
$\sigma_1 \ge \sigma_2 \ge \cdots \ge \sigma_N \ge 0$, and $U$ and $V$ are orthogonal matrices containing, respectively, the left and right singular vectors. The Moore--Penrose pseudoinverse of $A$ is then given by 
\begin{equation}
A^+ = V\,\Sigma^{-1}U^T\,.
\label{eq:pseudoinv}
\end{equation}
When the smallest singular values $\sigma_i$ approach zero, their reciprocals $1/\sigma_i$ grow without bound. As a result, numerical errors and noise in the data are greatly amplified in the computed solution. The modes associated with small singular values typically correspond to oscillatory, high-frequency components; thus, even if the exact solution is smooth, the discrete approximation may exhibit spurious oscillations and instability.

\subsection{Stabilizing the Solution: Truncated SVD and Tikhonov Regularization}
\label{sec:stabilization}
To mitigate the deleterious effects of small singular values, regularization is essential \cite{Hansen1998,Tikhonov1977}. Two widely used strategies are truncated singular value decomposition (TSVD) and Tikhonov regularization.

\textbf{Truncated SVD (TSVD).} In the TSVD approach, only the largest $r \le N$ singular values are retained, and the solution is approximated as 
\begin{equation}
\mathbf{\alpha}_r = \sum_{i=1}^r \frac{1}{\sigma_i}\, v_i\, (u_i^T \mathbf{f})\,,
\label{eq:tsvd_solution}
\end{equation}
where $u_i$ and $v_i$ are the $i$th columns of $U$ and $V$, respectively. This procedure acts as a low-pass filter: it suppresses components corresponding to small singular values (which typically encode high-frequency noise), while retaining the dominant modes that capture the smooth structure of the solution.

\textbf{Tikhonov regularization.} An alternative approach is to solve the penalized least-squares problem
\begin{equation}
\min_{\mathbf{\alpha}} \;\;\|A\mathbf{\alpha} - \mathbf{f}\|_2^2 + \lambda^2 \,\|\mathbf{\alpha}\|_2^2\,,
\label{eq:tikhonov_problem}
\end{equation}
where $\lambda > 0$ is the regularization parameter controlling the trade-off between fidelity to the data and smoothness of the solution. The normal equations for \eqref{eq:tikhonov_problem} are 
\begin{equation}
(A^T A + \lambda^2 I)\,\mathbf{\alpha} = A^T \mathbf{f}\,,
\label{eq:normal_eq}
\end{equation}
and the Tikhonov solution in SVD form is 
\begin{equation}
\mathbf{\alpha}_\lambda = \sum_{i=1}^{N} \frac{\sigma_i}{\sigma_i^2 + \lambda^2}\,(u_i^T \mathbf{f})\, v_i\,,
\label{eq:tikh_solution}
\end{equation}
which reveals the filter factors $\sigma_i/(\sigma_i^2 + \lambda^2)$ that continuously attenuate the influence of each mode (in contrast to the all-or-nothing truncation of TSVD).

Once the regularized coefficient vector $\mathbf{\alpha} = (\alpha^{(1)},\ldots,\alpha^{(N)})^T$ is obtained (either via TSVD or Tikhonov), the approximate solution of the Helmholtz equation can be reconstructed as 
\begin{equation}
u(X) \approx \sum_{j=1}^N \alpha^{(j)}\, \phi(\|X - X_j^c\|_2,\,\varepsilon)\,,
\label{eq:u_reconstruction}
\end{equation}
using the basis expansion \eqref{eq:u_expansion}. This MQ-RBF framework provides a high-order, meshfree method for solving the Helmholtz equation, with ill-posedness addressed through robust regularization techniques rooted in numerical linear algebra. The performance of the method depends critically on three factors: the choice of shape parameter $\varepsilon$, the conditioning of the RBF matrices, and the proper selection of the regularization parameter $\lambda$. These aspects will be central to the numerical experiments discussed later in the paper.

\section{Krylov Subspace and Tikhonov Regularization Methods}
\label{sec:krylov}
The SVD analysis in the previous section shows that small singular values of the MQ-RBF discretization matrix amplify noise and modeling errors. In practice, effective methods for MQ systems therefore blend regularization (to control the influence of small singular values) with iterative Krylov methods (to reduce computational cost while preserving stability). We summarize the ingredients used in our approach.

\subsection{Inexpensive Truncated SVD (Ine-TSVD)}
In practice, Ine-TSVD replaces the expensive full SVD of $A$ with a projected SVD built by a short Golub–Kahan bidiagonalization (GKB) \cite{Golub1965} and related Lanczos-based implementations (see also \cite{Hansen1998,Bjork1996}). After $\ell$ steps of Algorithm~\ref{alg:GKB}, we obtain orthonormal bases
$W_{\ell+1}\in\mathbb{R}^{N\times(\ell+1)}$ and $Z_\ell\in\mathbb{R}^{N\times\ell}$
together with a bidiagonal surrogate $C_\ell\in\mathbb{R}^{(\ell+1)\times\ell}$ satisfying
\begin{equation*}
AZ_\ell \;=\; W_{\ell+1}C_\ell,
\qquad
A^{\top}W_{\ell+1} \;=\; Z_\ell C_\ell^{\top}.
\end{equation*}
The singular values of $C_\ell$ (Ritz singular values of $A$) approximate the largest
singular values of $A$. Computing the small SVD $C_\ell=\widetilde U\,\widetilde\Sigma\,\widetilde V^{\top}$
yields approximate singular triplets of $A$ as
\begin{equation*}
\widehat u_i \;=\; W_{\ell+1}\widetilde u_i,\qquad
\widehat v_i \;=\; Z_\ell\widetilde v_i,\qquad
\widehat\sigma_i \;=\; \tilde\sigma_i,\qquad i=1,\dots,\ell.
\end{equation*}
Restricting to the first $k$ Ritz modes gives the projected TSVD solution
\begin{equation*}
\boldsymbol{\alpha}_k
\;=\;
\sum_{i=1}^{k} \frac{\widehat u_i^{\top} f}{\widehat\sigma_i}\,\widehat v_i
\;=\;
Z_\ell\,\widetilde V_k\,\widetilde\Sigma_k^{-1}\,\widetilde U_k^{\top}\,(\beta_1 e_1),
\end{equation*}
where $\beta_1=\|f\|_2$ and $e_1=(1,0,\dots,0)^{\top}\in\mathbb{R}^{\ell+1}$. In practice,
$\ell$ is taken modest (typically $\ell\!\approx\!2k$--$3k$), and $k$ is chosen by a simple
spectral criterion on the projected problem (e.g., Picard plot, visible spectral gap, or a
discrepancy rule on $\widehat\sigma_i$). Each GKB step requires one multiply with $A$ and one
with $A^{\top}$; the total work is
$\mathcal{O}\!\big(\ell\,\mathrm{mv}(A)+\ell\,\mathrm{mv}(A^{\top})+\ell^3\big)$
(where $\mathrm{mv}(\cdot)$ denotes the cost of a matrix--vector product),
with storage $\mathcal{O}(2N\ell)$. Optional (partial) reorthogonalization preserves the
numerical orthogonality of $W_{\ell+1}$ and $Z_\ell$; the process naturally halts on (near)
breakdown when some $\alpha_j$ or $\beta_{j+1}$ vanishes. In summary, Ine-TSVD reproduces the
low-pass filtering of TSVD at the cost of a short Krylov run and a tiny dense SVD, which is
particularly attractive for the dense MQ matrices produced by our RBF discretization.

\begin{algorithm}[hbt!]
\caption{Golub--Kahan bidiagonalization (Lanczos-based)}
\label{alg:GKB}
\begin{algorithmic}[1]
\Require $A\in\mathbb{R}^{N\times N}$, right-hand side $f\in\mathbb{R}^N$, number of steps $\ell$
\State $\beta_1 \gets \|f\|_2$, $\; w_1 \gets f/\beta_1$; set $z_0=\mathbf{0}$, $\beta_0=0$
\For{$j=1,2,\dots,\ell$}
  \State $r \gets A^T w_j - \beta_j z_{j-1}$ 
  \State $\alpha_j \gets \|r\|_2$;\quad $z_j \gets r/\alpha_j$
  \State $p \gets A z_j - \alpha_j w_j$
  \State $\beta_{j+1} \gets \|p\|_2$;\quad $w_{j+1} \gets p/\beta_{j+1}$
  \State (optional) reorthogonalize $w_{j+1}$ against $w_1,\dots,w_j$ and $z_j$ against $z_1,\dots,z_{j-1}$
\EndFor
\State Form $W_{\ell+1}=[w_1,\dots,w_{\ell+1}]\in\mathbb{R}^{N\times(\ell+1)}$, $Z_\ell=[z_1,\dots,z_\ell]\in\mathbb{R}^{N\times \ell}$,
and the lower bidiagonal $C_\ell\in\mathbb{R}^{(\ell+1)\times \ell}$ with diagonal $(\alpha_1,\dots,\alpha_\ell)$ and subdiagonal $(\beta_2,\dots,\beta_{\ell+1})$.
\State The relations $AZ_\ell = W_{\ell+1} C_\ell$ and $A^T W_{\ell+1} = Z_\ell C_\ell^T$ hold.
\end{algorithmic}
\end{algorithm}

\begin{algorithm}[hbt!]
\caption{Inexpensive TSVD (Ine--TSVD) using GKB}
\label{alg:IneTSVD}
\begin{algorithmic}[1]
\Require $A\in\mathbb{R}^{N\times N}$, right-hand side $f\in\mathbb{R}^N$, GKB steps $\ell$ (ex. $\ell\!\approx\!2k$--$3k$), truncation rank $k\le \ell$
\State Run Algorithm~\ref{alg:GKB} with input $(A,f,\ell)$ to obtain $W_{\ell+1}\in\mathbb{R}^{N\times(\ell+1)}$, $Z_\ell\in\mathbb{R}^{N\times \ell}$, $C_\ell\in\mathbb{R}^{(\ell+1)\times \ell}$, and $\beta_1=\|f\|_2$ with $w_1=f/\beta_1$
\State Compute the thin SVD of $C_\ell$: \quad $C_\ell=\widetilde U\,\widetilde \Sigma\,\widetilde V^T$, \quad
$\widetilde \Sigma=\mathrm{diag}(\tilde\sigma_1\ge\dots\ge \tilde\sigma_\ell)$
\State (Ritz triplets) For $i=1,\dots,\ell$, set \quad
$\widehat u_i \gets W_{\ell+1}\,\widetilde u_i$, \quad $\widehat v_i \gets Z_\ell\,\widetilde v_i$, \quad $\widehat \sigma_i \gets \tilde\sigma_i$
\State (Rank selection) Choose $k\in\{1,\dots,\ell\}$ (ex. via Picard plot, discrepancy, ou GCV sur les $\widehat \sigma_i$)
\State (TSVD on subspace) Using $f=\beta_1 w_1$ (donc $W_{\ell+1}^T f=\beta_1 e_1$), compute the reduced TSVD solution
\[
y_k \;\gets\; \sum_{i=1}^{k} \frac{\widetilde u_i^T(\beta_1 e_1)}{\tilde\sigma_i}\,\widetilde v_i
\;=\; \widetilde V_k\,\widetilde\Sigma_k^{-1}\,\widetilde U_k^T\,(\beta_1 e_1),
\]
then lift to the full space \quad $\,\alpha_k \gets Z_\ell\,y_k \;=\; \sum_{i=1}^{k} \frac{\widehat u_i^T f}{\widehat\sigma_i}\,\widehat v_i$.
\State \Return $\alpha_k$.
\end{algorithmic}
\end{algorithm}
\subsection{Tikhonov Regularization and Parameter Choice}
Instead of directly solving the ill-posed least-squares problem 
\(\min_{\boldsymbol{\alpha}} \|A\boldsymbol{\alpha}-\mathbf{f}\|_2^2\),
classical Tikhonov regularization \cite{Tikhonov1977} stabilizes the solution via
\begin{equation}
\min_{\boldsymbol{\alpha}}\ \|A\boldsymbol{\alpha}-\mathbf{f}\|_2^2+\lambda^2\|\boldsymbol{\alpha}\|_2^2,
\label{eq:gen_tikhonov}
\end{equation}
where \(\lambda>0\) balances fidelity and stability. The optimality condition yields the (well-conditioned) normal equations
\begin{equation}
(A^{\top}A+\lambda^2 I)\,\boldsymbol{\alpha}=A^{\top}\mathbf{f}.
\end{equation}
Let \(A=U\Sigma V^{\top}\) be an SVD with singular values \(\sigma_1\ge\dots\ge\sigma_r>0\) (\(r=\mathrm{rank}(A)\)), and define
\(\widehat{\mathbf f}=U^{\top}\mathbf f\) and the Tikhonov filter factors
\begin{equation}
\phi_i(\lambda)=\frac{\sigma_i^2}{\sigma_i^2+\lambda^2}\in(0,1),\qquad i=1,\dots,r.
\end{equation}
Then the Tikhonov solution is
\begin{equation}
\boldsymbol{\alpha}_{\lambda}
=\sum_{i=1}^{r}\phi_i(\lambda)\,\frac{\widehat f_i}{\sigma_i}\,v_i,
\label{eq:tikh_svd_solution}
\end{equation}
which shows that small-\(\sigma_i\) modes are smoothly damped (contrast with the hard truncation of TSVD).

\paragraph{Residual and solution norms in SVD variables}
With \(\phi_i(\lambda)\) as above,
\begin{align}
\rho^2(\lambda)&:=\|A\boldsymbol{\alpha}_\lambda-\mathbf f\|_2^2
= \sum_{i=1}^{r}\bigl(1-\phi_i(\lambda)\bigr)^2\,\widehat f_i^{\,2}
+\sum_{i=r+1}^{N}\widehat f_i^{\,2}, \label{eq:rho}\\
\eta^2(\lambda)&:=\|\boldsymbol{\alpha}_\lambda\|_2^2
= \sum_{i=1}^{r}\left(\frac{\phi_i(\lambda)}{\sigma_i}\right)^2 \widehat f_i^{\,2}.
\label{eq:eta}
\end{align}

\paragraph{Generalized cross-validation (GCV)}
Define the matrix 
\[
H(\lambda) = A\,(A^{\top}A + \lambda^{2}I)^{-1}A^{\top},
\]
whose trace can be expressed as 
\(\mathrm{tr}\,H(\lambda) = \sum_{i=1}^{r}\phi_i(\lambda)\).
A numerically stable formulation of the generalized cross–validation (GCV) \cite{Golub1979,Wahba1990,Hansen1998} function is given by :
\begin{equation}
\mathrm{GCV}(\lambda)=
\frac{\|(I-H(\lambda))\mathbf f\|_2^2}{\bigl[N-\mathrm{tr}\,H(\lambda)\bigr]^2}
=\frac{\displaystyle\sum_{i=1}^{r}\bigl(1-\phi_i(\lambda)\bigr)^2\,\widehat f_i^{\,2}
+\sum_{i=r+1}^{N}\widehat f_i^{\,2}}
{\displaystyle\biggl(N-\sum_{i=1}^{r}\phi_i(\lambda)\biggr)^2}.
\label{eq:gcv_svd}
\end{equation}
We select \(\lambda\) that minimizes \(\mathrm{GCV}(\lambda)\). This requires only \(\{\sigma_i,\widehat f_i\}\), not explicit formation of \(A^{\top}A\).

\paragraph{L-curve criterion (detailed)}
The L-curve \cite{Hansen1992,Hansen1993,Hansen1998} is the parametric plot
\[
\mathcal L(\lambda)=\bigl(x(\lambda),y(\lambda)\bigr):=
\bigl(\log \rho(\lambda),\ \log \eta(\lambda)\bigr),
\]
with \(\rho,\eta\) given by \eqref{eq:rho}–\eqref{eq:eta}. The recommended \(\lambda\) is the
\emph{corner} of \(\mathcal L\), where curvature is maximal, balancing residual fit and solution smoothness.

\begin{algorithm}[hbt!]
\caption{Classical Tikhonov with parameter (Tikh-Reg) via GCV}
\label{alg:TikhRgGCV}
\begin{algorithmic}[1]
\Require $A\in\mathbb{R}^{N\times N}$, data $\mathbf f\in\mathbb{R}^N$, log-grid $\Lambda=\{\lambda_j\}_{j=1}^J$
\State Compute economy SVD: $A=U\Sigma V^{\top}$ with $\Sigma=\mathrm{diag}(\sigma_1,\dots,\sigma_r)$
\State $\widehat{\mathbf f}=U^{\top}\mathbf f$, \quad $s=\|\mathbf f\|_2^2$, \quad $s_r=\sum_{i=1}^{r}\widehat f_i^{\,2}$, \quad $s_\perp=s-s_r$
\For{$j=1,\dots,J$} \Comment{evaluate GCV on the SVD}
  \State $\lambda\leftarrow\lambda_j$, \quad $\phi_i(\lambda)=\dfrac{\sigma_i^2}{\sigma_i^2+\lambda^2}$ for $i=1,\dots,r$
  \State Residual: $\displaystyle \rho^2(\lambda_j)=\sum_{i=1}^{r}\bigl(1-\phi_i(\lambda)\bigr)^2\,\widehat f_i^{\,2}\;+\;s_\perp$
  \State Trace term: $\displaystyle \tau(\lambda_j)=\sum_{i=1}^{r}\phi_i(\lambda)$
  \State $\displaystyle \mathrm{GCV}(\lambda_j)=\frac{\rho^2(\lambda_j)}{\bigl[N-\tau(\lambda_j)\bigr]^2}$
\EndFor
\State $\lambda_\ast=\arg\min_{\lambda_j\in\Lambda}\ \mathrm{GCV}(\lambda_j)$
\State Solution at $\lambda_\ast$: $\displaystyle \boldsymbol{\alpha}_{\lambda_\ast}
=\sum_{i=1}^{r}\frac{\phi_i(\lambda_\ast)}{\sigma_i}\,\widehat f_i\,v_i
\;=\;V\,\mathrm{diag}\!\Big(\frac{\phi_1(\lambda_\ast)}{\sigma_1},\dots,\frac{\phi_r(\lambda_\ast)}{\sigma_r}\Big)\,\widehat{\mathbf f}_{1:r}$
\State \Return $\boldsymbol{\alpha}_{\lambda_\ast}$, $\lambda_\ast$.
\end{algorithmic}
\end{algorithm}

\subsection{Hybrid Krylov--Tikhonov via Golub--Kahan Projection}
Golub--Kahan bidiagonalization (GKB)  produces orthonormal bases
\(W_{\ell+1}\in\mathbb{R}^{N\times(\ell+1)}\), \(Z_\ell\in\mathbb{R}^{N\times\ell}\),
and a lower-bidiagonal matrix \(C_\ell\in\mathbb{R}^{(\ell+1)\times\ell}\) such that
\begin{equation}
A Z_\ell = W_{\ell+1}\, C_\ell,
\qquad
A^{\top} W_{\ell+1} = Z_\ell\, C_\ell^{\top}.
\label{eq:GKB}
\end{equation}
For modest \(\ell\), \(\mathrm{range}(W_{\ell+1})\) and \(\mathrm{range}(Z_\ell)\) capture the dominant SVD subspaces of \(A\) to high accuracy (see section 3.1). 
Projecting Tikhonov onto \(\mathrm{range}(Z_\ell)\) yields the hybrid problem 
\begin{equation}
\min_{\mathbf y\in\mathbb{R}^{\ell}}\ \|C_\ell \mathbf y - W_{\ell+1}^{\top}\mathbf f\|_2^2
+\lambda^2 \|\mathbf y\|_2^2,
\qquad \boldsymbol{\alpha}=Z_\ell \mathbf y.
\label{eq:hybrid_problem}
\end{equation}

\paragraph{Cheap GCV/L-curve on the projected problem}
Because \eqref{eq:hybrid_problem} lives in \(\mathbb{R}^{\ell}\) with data in \(\mathbb{R}^{\ell+1}\),
parameter-choice rules can be evaluated on a \emph{tiny} surrogate, making them far less expensive than on \(A\).
Writing \(W_{\ell+1}^{\top}\mathbf f=\beta_1 e_1\) (as in Algorithm~\ref{alg:GKB}),
let \(C_\ell=\widetilde U\,\widetilde\Sigma\,\widetilde V^{\top}\) be the SVD with singular values
\(\tilde\sigma_1\ge\dots\ge \tilde\sigma_\ell\ge 0\).
Define projected filter factors \(\tilde\phi_i(\lambda)=\tilde\sigma_i^2/(\tilde\sigma_i^2+\lambda^2)\) and
\(\widetilde{\mathbf f}=\widetilde U^{\top}(\beta_1 e_1)\).
Then the residual and solution norms for the projected problem are
\begin{align}
\rho_\ell^2(\lambda)&=\|C_\ell \mathbf y_\lambda - \beta_1 e_1\|_2^2
=\sum_{i=1}^{\ell}\bigl(1-\tilde\phi_i(\lambda)\bigr)^2\,\widetilde f_i^{\,2}
+\widetilde f_{\ell+1}^{\,2},\label{eq:rho_proj}\\
\eta_\ell^2(\lambda)&=\|\mathbf y_\lambda\|_2^2
=\sum_{i=1}^{\ell}\left(\frac{\tilde\phi_i(\lambda)}{\tilde\sigma_i}\right)^2\,\widetilde f_i^{\,2}.
\label{eq:eta_proj}
\end{align}
\emph{GCV on the surrogate} uses the hat matrix
\(H_\ell(\lambda)=C_\ell(C_\ell^{\top}C_\ell+\lambda^2 I)^{-1}C_\ell^{\top}\) with
\(\mathrm{tr}\,H_\ell(\lambda)=\sum_{i=1}^{\ell}\tilde\phi_i(\lambda)\) and sets
\begin{equation}
\mathrm{GCV}_\ell(\lambda)
=\frac{\|(I-H_\ell(\lambda))(\beta_1 e_1)\|_2^2}{\bigl[(\ell+1)-\mathrm{tr}\,H_\ell(\lambda)\bigr]^2}
=\frac{\rho_\ell^2(\lambda)}{\bigl[(\ell+1)-\sum_{i=1}^{\ell}\tilde\phi_i(\lambda)\bigr]^2}.
\label{eq:gcv_proj}
\end{equation}
\emph{L-curve on the surrogate} is the parametric curve
\(\mathcal L_\ell(\lambda)=(\log\rho_\ell(\lambda),\log\eta_\ell(\lambda))\) from
\eqref{eq:rho_proj}–\eqref{eq:eta_proj}; the recommended \(\lambda_\ast\) is found at the
corner, e.g., by maximizing the discrete curvature on a logarithmic grid of \(\lambda\)
(as detailed in the previous subsection).

\begin{algorithm}[t]
\caption{Hybrid Krylov--Tikhonov(HKT) via GKB (parameter choice by GCV or L-curve)}
\label{alg:HybridTik}
\begin{algorithmic}[1]
\Require System \(A\in\mathbb{R}^{N\times N}\), data \(\mathbf f\in\mathbb{R}^N\), GKB steps \(\ell\),
selection rule \(\in\{\mathrm{GCV},\mathrm{L\text{-}curve}\}\)
\State Run Algorithm~\ref{alg:GKB} with \((H,\mathbf f,\ell)\) to get \(W_{\ell+1},Z_\ell,C_\ell\) and \(\beta_1=\|\mathbf f\|_2\) with \(W_{\ell+1}^{\top}\mathbf f=\beta_1 e_1\)
\State Compute SVD \(C_\ell=\widetilde U\,\widetilde\Sigma\,\widetilde V^{\top}\); set \(\widetilde{\mathbf f}=\widetilde U^{\top}(\beta_1 e_1)\)
\State Define \(\tilde\phi_i(\lambda)=\tilde\sigma_i^2/(\tilde\sigma_i^2+\lambda^2)\), and compute
\(\rho_\ell(\lambda),\eta_\ell(\lambda)\) from \eqref{eq:rho_proj}–\eqref{eq:eta_proj}
\If{selection rule = GCV}
  \State Choose \(\lambda_\ast=\arg\min_{\lambda>0}\ \mathrm{GCV}_\ell(\lambda)\) using \eqref{eq:gcv_proj}
\Else
  \State Choose \(\lambda_\ast\) at the L-curve corner (max curvature) computed from \(\{\rho_\ell(\lambda),\eta_\ell(\lambda)\}\) on a log-grid
\EndIf
\State Solve \((C_\ell^{\top}C_\ell+\lambda_\ast^2 I)\,\mathbf y_{\lambda_\ast}=C_\ell^{\top}(\beta_1 e_1)\)
\State Set \(\boldsymbol{\alpha}=Z_\ell \mathbf y_{\lambda_\ast}\) and \Return \(\boldsymbol{\alpha}\).
\end{algorithmic}
\end{algorithm}

\paragraph{Cost and storage}
The total work is \(O\!\big(\ell\,(\mathrm{mv}(A)+\mathrm{mv}(A^{\top}))+\ell^3\big)\); storage is
\(O(2N\ell)\).
Evaluating GCV or the L-curve on the surrogate replaces full-scale scans by operations on
\(\ell\times\ell\) matrices and vectors of length \(\ell+1\), yielding substantial savings for
dense MQ matrices.

\subsection{Krylov Methods as Iterative Regularization: LSQR and GMRES}
Krylov methods can themselves act as implicit regularizers. LSQR \cite{Paige1982} is based on GKB and is algebraically equivalent to conjugate gradients on the normal equations, yet avoids explicitly forming $A^T A$. Let $\mathbf{\alpha}_k$ denote the $k$th LSQR iterate. Because Krylov spaces emphasize the dominant right singular subspace first, early iterates $\mathbf{\alpha}_k$ primarily combine large-$\sigma_i$ modes (which are stable), while small-$\sigma_i$ modes (which are noise-amplifying) enter only at later iterations. This produces the well-known \emph{semi-convergence} behavior: the error decreases up to an optimal iteration $k_{\text{opt}}$ and then increases as noise-dominated modes start to contaminate the solution \cite{Hansen1998}. A practical regularization strategy is therefore \emph{early stopping} at $k_{\text{opt}}$, chosen by, e.g., the L-curve criterion applied to the sequence $\{ (\|A \mathbf{\alpha}_k - \mathbf{f}\|_2,\; \|\mathbf{\alpha}_k\|_2) \}$ or GCV variants adapted to Krylov filtering (see \cite{Hansen1998,Hansen1992}). For nonsymmetric formulations one may also employ GMRES \cite{Saad1986} (possibly with right preconditioning), which exhibits a similar filtering effect on the spectrum of $A$.
\section{Existence and Convergence Analysis}
\label{sec:existence-convergence}
This section establishes (i) invertibility statements for multiquadric matrices used in interpolation and collocation; (ii) existence and uniqueness of regularized discrete solutions; and (iii) convergence (and smoothing) of the proposed methods.

\subsection{Invertibility of MQ matrices}

We recall that the MQ kernel $\phi(r,\varepsilon)=\sqrt{1+\varepsilon^2 r^2}$ is \emph{conditionally positive definite of order~$1$} (CPD(1)), and therefore its distance matrix is invertible on the subspace orthogonal to constants (or, equivalently, after the standard polynomial augmentation of degree~$0$). See Micchelli's seminal result \cite{Micchelli1986} and the expositions in \cite[Ch.~2]{Sarra2009}.  A concise formulation adapted to MQ is:

\begin{theorem}[Micchelli's CPD(1) criterion specialized to MQ] \cite{Micchelli1986}
\label{thm:micchelli-mq}
Let $\{X_j\}_{j=1}^N\subset\mathbb{R}^d$ be pairwise distinct. Define the MQ distance matrix $B\in\mathbb{R}^{N\times N}$ by $B_{ij}=\phi(\|X_i-X_j\|_2,\varepsilon)$. Then $B$ is \emph{conditionally} positive definite of order $1$, i.e.,
\[
\mathbf{c}^T B \mathbf{c} > 0 \quad\text{for all }\mathbf{c}\neq 0\text{ with }\sum_{j=1}^N c_j=0.
\]
Consequently, the augmented system
\[
\begin{bmatrix} B & \mathbf{1} \\ \mathbf{1}^T & 0 \end{bmatrix}
\begin{bmatrix} \boldsymbol{\alpha} \\ \lambda \end{bmatrix}
=
\begin{bmatrix} \mathbf{f} \\ 0 \end{bmatrix}
\]
is nonsingular and yields a unique interpolant.
\end{theorem}
\emph{Proof.} This is a direct application of Micchelli's theorem for CPD kernels with order $m=1$. For MQ one checks that $\psi(r):=\phi(\sqrt{r},\varepsilon)$ has completely monotone derivative on $(0,\infty)$, and thus the CPD(1) property holds; see the derivative test and worked-out derivatives in \cite[Sec.~2.3]{Sarra2009}. \qed
\begin{lemma}[Complete monotonicity for MQ]
\label{lem:cm-mq}
Let $\psi(r)=\phi(\sqrt{r},\varepsilon)=\sqrt{1+\varepsilon^2 r}$. Then $(-1)^\ell \psi^{(\ell)}(r)\ge 0$ for all integers $\ell\ge 1$ and $r>0$. \\\emph{Proof.} Differentiate explicitly:
\[
\psi'(r)=\frac{\varepsilon^2}{2(1+\varepsilon^2 r)^{1/2}},\quad
\psi''(r)=-\frac{\varepsilon^4}{4(1+\varepsilon^2 r)^{3/2}},\quad
\psi^{(3)}(r)=\frac{3\varepsilon^6}{8(1+\varepsilon^2 r)^{5/2}},\ \ldots
\]
which alternate in sign, establishing complete monotonicity of $\psi'$, hence CPD(1) \\for MQ \cite[Sec.~2.3]{Sarra2009}.\qed
\end{lemma}

\begin{proposition}[Generalized MQ and polynomial augmentation]
\label{prop:gmq-augmentation}
For the generalized MQ $\phi_\beta(r;\varepsilon)=(1+\varepsilon^2 r^2)^\beta$ with $\beta\notin\mathbb{N}_0$, the kernel is strictly positive definite when $\beta<0$ and CPD of order $\lceil \beta\rceil$ when $\beta>0$. In the latter case, nonsingularity of the interpolation system is guaranteed after augmenting with all polynomials of degree $< \lceil \beta\rceil$ (standard side constraints). \\\emph{Proof.} See \cite[Sec.~2.10]{Sarra2009} and \cite{Micchelli1986,Madych1990}. \qed
\end{proposition}

\begin{remark}[Ill-conditioning vs.~accuracy]
Even though the collocation matrix is invertible under the hypotheses above, the conditioning deteriorates as $\varepsilon\to 0$ (flat basis) or as the fill distance decreases. This is the well-known RBF ``uncertainty principle'': accuracy tends to improve as conditioning worsens \cite[Sec.~2.5]{Sarra2009}. \qed
\end{remark}

\subsection{Discrete collocation for Helmholtz and existence with regularization}

Let $A\in\mathbb{R}^{N\times N}$ denote the (unsymmetric) Kansa collocation matrix enforcing the PDE at interior points and boundary data at boundary points (cf.~Section~\ref{sec:MQ}). For MQ bases, $A$ is dense and often severely ill-conditioned. In this paper we stabilize the discrete problems using only three key algorithms: (i) Inexpensive Truncated SVD  via a short Golub--Kahan projection; (ii) Classical Tikhonov regularization with parameter choice by GCV or the L-curve; and (iii) Hybrid Krylov--Tikhonov via GKB, which evaluates GCV/L-curve on a small projected problem. Regardless of the rank of $A$, the Tikhonov and TSVD formulations below produce a \emph{unique} regularized solution.

\begin{proposition}[Existence and uniqueness for Tikhonov]
\label{prop:tikh-existence}
For any $\lambda>0$, the Tikhonov minimizer
\[
\min_{y\in\mathbb{R}^N}\ \|Ay-f\|_2^2+\lambda^2\|y\|_2^2
\]
exists, is unique, and satisfies $(A^{\top}A+\lambda^2 I)\,y_\lambda=A^{\top}f$.\\
\emph{Proof.} $A^{\top}A+\lambda^2 I$ is symmetric positive definite for $\lambda>0$. \qed
\end{proposition}

\begin{proposition}[Existence and uniqueness for TSVD]
\label{prop:tsvd-existence}
Let $A=U\Sigma V^{\top}$ be an SVD and fix $k\in\{1,\ldots,\mathrm{rank}(A)\}$. The TSVD solution
\(
y_k=\sum_{i=1}^k \sigma_i^{-1} (u_i^{\top} f)\, v_i
\)
is the unique minimum-norm solution in $\mathrm{range}(V_k)$ to
$\min_{y\in\mathrm{range}(V_k)}\|A\alpha-f\|_2$.\\
\emph{Proof.} Standard spectral projection argument \cite{Bjork1996,Hansen1998}. \qed
\end{proposition}

\begin{remark}[Projected (hybrid) formulations and cheap parameter choice]
After $\ell$ steps of GKB started with $f$, we obtain
\[
A Z_\ell=W_{\ell+1}\,\overline{C}_\ell,\qquad
A^{\top} W_{\ell+1}=Z_\ell\,\overline{C}_\ell^{\top},\qquad W_{\ell+1}^{\top}f=\beta_1 e_1,
\]
with orthonormal $W_{\ell+1}\in\mathbb{R}^{N\times(\ell+1)}$, $Z_\ell\in\mathbb{R}^{N\times\ell}$ and a small bidiagonal
$\overline{C}_\ell\in\mathbb{R}^{(\ell+1)\times\ell}$. Solving Tikhonov/TSVD on the \emph{projected} problem
\[
\min_{y\in\mathbb{R}^{\ell}}\ \|\overline{C}_\ell y-\beta_1 e_1\|_2^2+\lambda^2\|y\|_2^2,\qquad x=Z_\ell y,
\]
is well-posed for any $\lambda>0$, inherits the dominant right singular subspace of $A$ captured by $Z_\ell$, and crucially allows GCV or L-curve to be evaluated on an $\ell\times\ell$ surrogate at negligible cost compared to the full matrix \cite{Bjork1996,Hansen1998,Lewis2009}. \qed
\end{remark}

\subsection{Convergence and smoothing properties}

Assume the standard \emph{discrete Picard condition} for discrete first-kind problems: the Fourier coefficients of the noiseless data in the left singular basis of $A$ decay faster (on average) than the singular values.

\begin{assumption}[Discrete Picard condition]
\label{ass:picard}
Let $A=U\Sigma V^{\top}$ and $f=f^\dagger+e$ with $\|e\|\le\delta$. There exist $\nu>0$ and $C>0$ such that
$|u_i^{\top} f^\dagger|\le C\,\sigma_i^{1+\nu}$ for all $i$.
\end{assumption}

\paragraph{Filter factors (only the three methods used)}
\[
\text{Tikhonov: }\ \varphi_i^{\rm Tik}(\lambda)=\frac{\sigma_i^2}{\sigma_i^2+\lambda^2},\qquad
\text{TSVD (Ine-TSVD): }\ \varphi_i^{\rm TSVD}(k)=\mathbf{1}_{\{i\le k\}},
\]
and, for the hybrid method, the same formulas apply with the \emph{Ritz} singular values of $\overline{C}_\ell$.

\begin{theorem}[Convergence of Tikhonov under Assumption~\ref{ass:picard}]
\label{thm:tikhonov-convergence}
Let $y_\lambda$ solve $\min_y\|Ay-f\|_2^2+\lambda^2\|y\|_2^2$.
If $\lambda=\lambda(\delta)\downarrow 0$ and $\delta/\lambda(\delta)\to 0$ as $\delta\downarrow 0$, then
$y_{\lambda(\delta)}\to y^\dagger$, the minimum-norm exact solution of $Ay=f^\dagger$
\cite{Hansen1998,Tikhonov1977}. Moreover, $y_\lambda$ is a \emph{smoothed} least-squares solution:
high-frequency (small-$\sigma_i$) components are damped by $\varphi_i^{\rm Tik}(\lambda)$ \cite{Hansen1998,Tikhonov1977}.
\end{theorem}

\begin{theorem}[Convergence of TSVD (Ine-TSVD)] 
\label{thm:tsvd-convergence}  Let $y_k$ be the TSVD solution. If $k=k(\delta)\to\infty$ and $\sigma_{k(\delta)}\gg\delta$ (e.g., $k$ by GCV),
then $y_{k(\delta)}\to y^\dagger$ as $\delta\to 0$ under Assumption.
The solution is \emph{smoothed} by the hard cut-off $\varphi_i^{\rm TSVD}(k)$ that removes noise-dominated modes \cite{Hansen1998,Golub1979,Bjork1996}.
\end{theorem}

\begin{theorem}[Convergence and smoothing of the Hybrid Krylov--Tikhonov (via GKB)]
\label{thm:hybrid}
Let $Z_\ell$ be the GKB right basis after $\ell$ steps, and suppose the dominant right singular subspace of $A$
up to index $k$ is contained in $\mathrm{range}(Z_\ell)$ (typically $\ell\approx 2k$--$3k$).
Let $\lambda$ be chosen on the projected problem by GCV or L-curve.
Then the hybrid solution $x_\ell=Z_\ell y_\lambda$ is \emph{quasi-optimal} among rank-$k$ spectral filters:
its error differs from the full Tikhonov error by the (small) subspace approximation error; the associated
filter factors act on the Ritz singular values of $\overline{C}_\ell$ \cite{Bjork1996,Lewis2009,Hansen1998}.
\end{theorem}
\section{Numerical Experiments}
\label{sec:experiments}

This section presents numerical experiments that evaluate the proposed Krylov–based methods against established approaches. In Examples~1 and~2 we solve the Helmholtz equation on a unit cube and a unit sphere, respectively, and compare the Hybrid Krylov–Tikhonov method (HKT) with classical Tikhonov regularization (Tikh-Rg) and a Tikhonov-regularized GMRES baseline (Reg-GMRES). Example~3 considers a complex industrial pump-casing geometry and benchmarks HKT against Reg-GMRES. Throughout all three examples, the regularization parameter is chosen \emph{exclusively} by generalized cross-validation (GCV).

Accuracy is reported using the relative error
\begin{equation}
\mathrm{Re} \;=\; \frac{\|u - u_a\|_\infty}{\|u\|_\infty},
\label{eq:rel_error}
\end{equation}
where \(u\) denotes the exact (manufactured) solution and \(u_a\) the numerical approximation; since \(u\) is available in every test, errors are computed directly. Computational efficiency is measured by elapsed CPU time (seconds).

All collocation points coincide with the nodes of three-dimensional finite element meshes generated using the open-source mesh generator \emph{Gmsh} \cite{GeuzaineRemacle2009}. For Examples~1 and~2 we consider \(N=359\), \(2154\), and \(6511\) collocation points obtained by successively refining the tetrahedral meshes on the cube and sphere. For Example~3 we use \(N=1190\) and \(6794\), corresponding to two refinement levels of the pump-casing mesh.\footnote{These values correspond to increasing refinement levels in the Gmsh-generated meshes used for the pump-casing geometry (Example~3).} The conditioning of the system matrix \(A\) (and thus the difficulty of the solve) increases with the wavenumber \(k\); to provide a consistent basis for comparison, we fix \(k=3\) in all runs.

All experiments are performed in \textsc{MATLAB} R2014a with IEEE double precision (machine epsilon \(\approx 10^{-16}\)) on a Windows~10 (64-bit) workstation equipped with an Intel\textsuperscript{\textregistered} Core\textsuperscript{TM} i5-3470 CPU @ 3.20~GHz and 8~GB RAM. For each method, we report relative error, runtime, and, when applicable, the GCV-selected parameter and iteration counts so that cost and accuracy can be compared on equal footing.
\subsection{Example 1: Helmholtz Equation in the Unit Cube}
\begin{figure}
\begin{center}
\includegraphics[width=.33\textwidth,height=.33\textwidth]{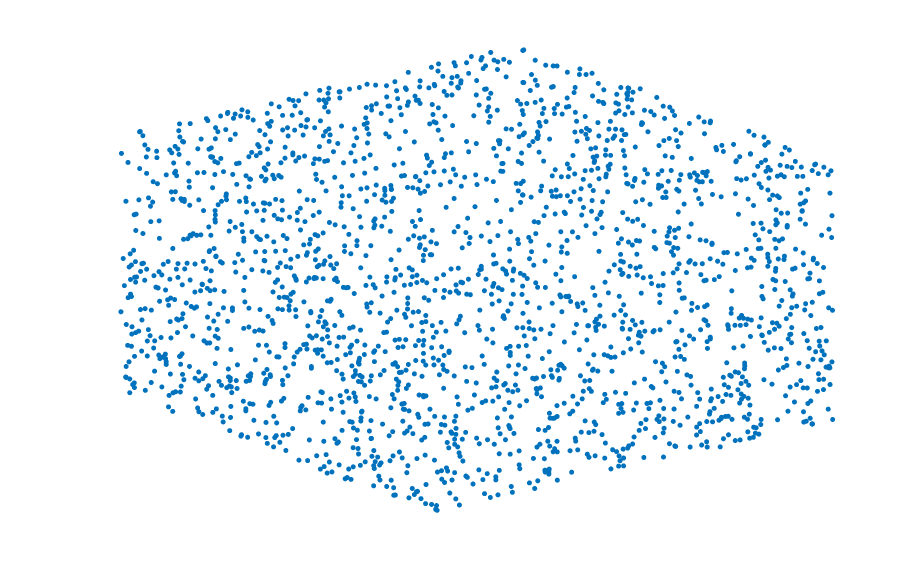}%
\includegraphics[width=.33\textwidth,height=.33\textwidth]{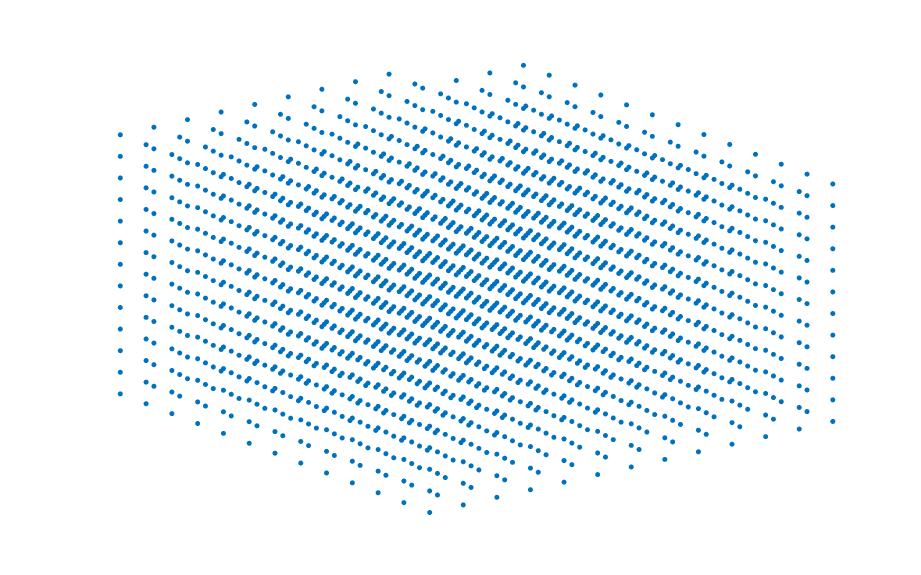}%
\includegraphics[width=.33\textwidth,height=.33\textwidth]{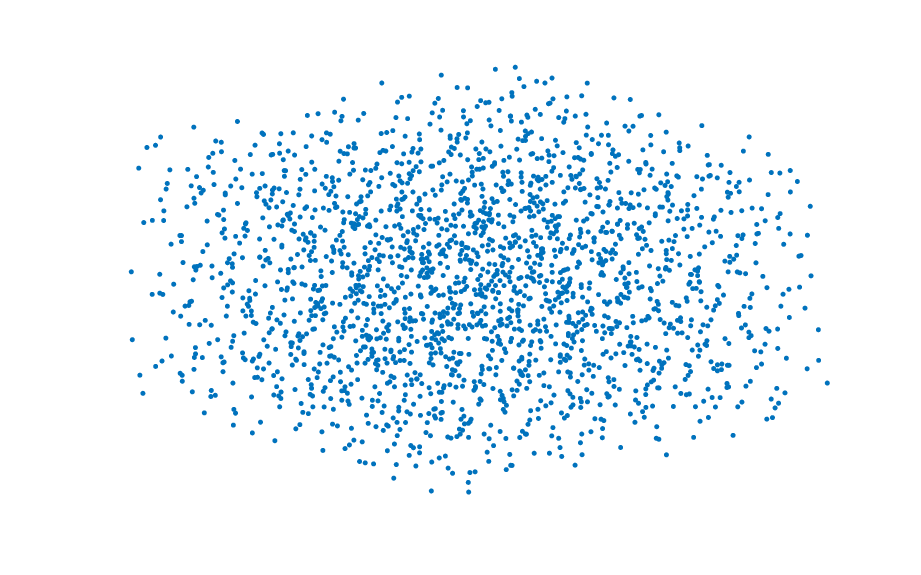}
\caption{Collocation points in the unit cube for Example~1: (left) random distribution, (middle) uniform grid distribution, (right) Halton sequence.}
\label{fig:points_cube}
\end{center}
\end{figure}

In this example, we consider the numerical solution of the three-dimensional Helmholtz problem \eqref{eq:Helmholtz} in the unit cube $\Omega = [0,1]^3$, illustrated in Figure~\ref{fig:points_cube}. As a test case, the source term $f$ is chosen such that the exact solution is 
\begin{equation}
u(x,y,z) = \exp\!\Big(-\frac{x^2 + y^2 + z^2}{\sigma}\Big)\,,
\qquad \sigma = 20\,,
\label{eq:exact_cube}
\end{equation}
and we impose Dirichlet boundary conditions using this exact solution on $\partial\Omega$. 

\begin{figure}
\begin{center}
\includegraphics[width=.5\textwidth,height=.4\textwidth]{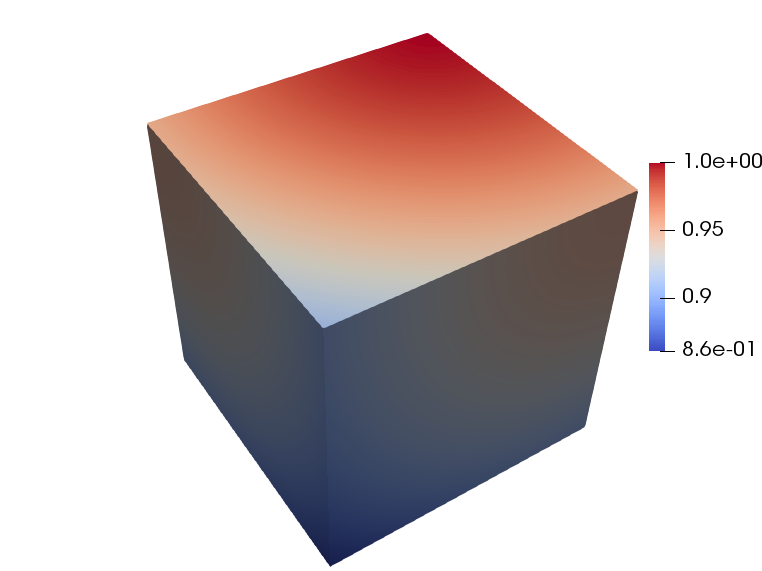}%
\includegraphics[width=.5\textwidth,height=.4\textwidth]{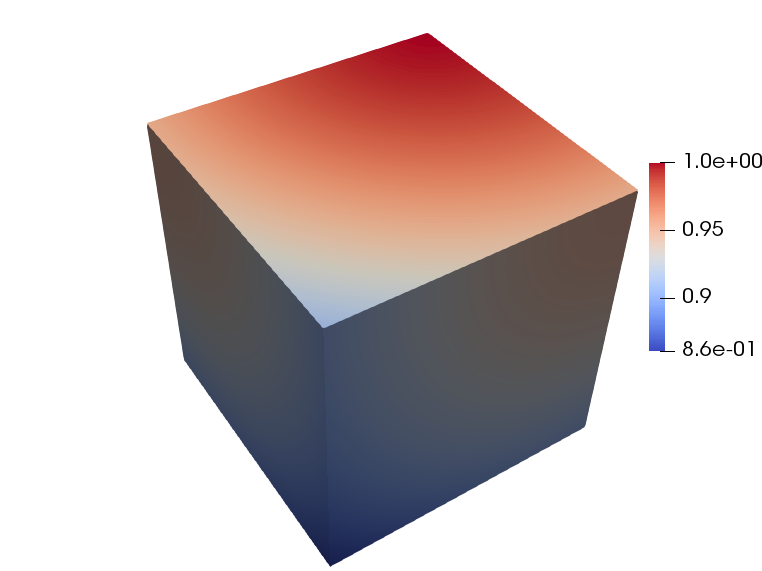}
\caption{Example~1 (Unit Cube): (Left) exact solution $u(x,y,z)$ given by \eqref{eq:exact_cube}, and (Right) the approximate solution obtained by the HKT(140) method for $k=3$ with $N=2154$ uniformly distributed points. The approximate solution visually matches the exact solution, indicating the effectiveness of the regularized Krylov method.}
    \label{fig:solution_cube}
\end{center}
\end{figure}

To assess the performance of different methods, we compare the inexpensive truncated SVD method (Ine-TSVD(140)), standard Tikhonov regularization (Tikh-Reg), and the Hybrid Krylov-Tikhonov method (HKT(140)) for different numbers $N$ of collocation points, with wavenumber $k = 3$. It is well known that the choice of node distribution affects the results of RBF-based meshless methods, so we consider \emph{random}, \emph{uniform}, and \emph{Halton} distributions of points. For comparison, we also include results obtained with Tikhonov-regularized GMRES(140) (abbreviated Reg-GMRES) \cite{Lewis2009,Saad1986,Calvetti2002}.
Figure~\ref{fig:points_cube} shows the collocation points in the cube for these distributions.
Figure~\ref{fig:solution_cube} shows a sample result of the exact solution and the approximate solution computed by the HKT(140) method for $k=3$ in the unit cube with $N=2154$ uniformly distributed points. Table~\ref{tab:cube_results} reports the efficiency and accuracy of the tested algorithms for various values of $N$. We list the relative errors (Re) and CPU times for each method under each point distribution. As expected, the errors generally decrease as $N$ increases (reflecting improved approximation with more points) and as the point distribution becomes more uniform (Halton sequences giving slightly better accuracy than random distributions). HKT consistently achieves highly accurate solutions with moderate cost.  
Ine-TSVD is generally faster but can lose accuracy for large $N$ if significant modes are truncated.  
Tikh-Reg produces smooth, stable solutions comparable to HKT but at higher computational cost due to the full SVD or normal equations.  
Reg-GMRES shows semi-convergence: errors initially decrease but increase again if iterations continue, confirming the theoretical discussion in Section~\ref{sec:existence-convergence}. 
When early stopping is applied at the optimal iteration, Reg-GMRES attains reasonable accuracy but remains less efficient than HKT.
\begin{table}[ht]
\centering
\caption{Unit cube, $k = 3$ (Example~1). Relative errors (Re) and CPU time (s). $\mathrm{Re}$ and $\kappa(A)$ are rounded.}
\label{tab:cube_results}
\tiny
\setlength{\tabcolsep}{5pt}
\renewcommand{\arraystretch}{1.15}
\rowcolors{4}{gray!4}{white}
\begin{tabular}{l c c cc cc cc cc}
\toprule
Point distribution & $N$ & $\kappa(A)$
& \multicolumn{2}{c}{Tikh-Rg} & \multicolumn{2}{c}{Ine-TSVD} & \multicolumn{2}{c}{HKT}
& \multicolumn{2}{c}{Reg-GMRES} \\
\cmidrule(lr){4-5}\cmidrule(lr){6-7}\cmidrule(lr){8-9}\cmidrule(lr){10-11}
 &  &  & Re & CPU & Re & CPU & Re & CPU & Re & CPU \\
\midrule
Random  & 359  & $2.9\times10^{18}$ & $3.1\times10^{-4}$ & 0.27 & $9.6\times10^{-5}$ & 0.18 & $1.4\times10^{-5}$ & 0.12 & $1.8\times10^{-1}$ & 0.29 \\
        & 2154 & $1.0\times10^{20}$ & $3.0\times10^{-5}$ & 7.78 & $5.7\times10^{-6}$ & 5.07 & $5.8\times10^{-6}$ & 4.66 & $1.6\times10^{-1}$ & 6.57 \\
        & 6511 & $2.8\times10^{21}$ & $2.2\times10^{-7}$ & 85.54& $1.2\times10^{-6}$ & 32.36& $1.2\times10^{-6}$ & 29.97& $1.3\times10^{-1}$ & 38.19 \\
\addlinespace[2pt]
Uniform & 359  & $7.1\times10^{17}$ & $2.7\times10^{-4}$ & 0.18 & $1.3\times10^{-3}$ & 0.15 & $3.0\times10^{-4}$ & 0.06 & $1.0\times10^{-1}$ & 0.36 \\
        & 2154 & $1.0\times10^{20}$ & $5.0\times10^{-5}$ & 7.51 & $4.8\times10^{-5}$ & 4.67 & $4.8\times10^{-5}$ & 4.47 & $3.0\times10^{-1}$ & 6.43 \\
        & 6511 & $3.5\times10^{20}$ & $6.4\times10^{-6}$ & 88.49& $1.6\times10^{-4}$ & 30.56& $1.5\times10^{-4}$ & 31.95& $5.9\times10^{-1}$ & 40.70 \\
\addlinespace[2pt]
Halton  & 359  & $2.9\times10^{18}$ & $4.5\times10^{-7}$ & 0.12 & $9.6\times10^{-5}$ & 0.14 & $8.1\times10^{-7}$ & 0.07 & $8.3\times10^{-2}$ & 0.24 \\
        & 2154 & $9.0\times10^{19}$ & $1.5\times10^{-7}$ & 7.23 & $5.3\times10^{-7}$ & 4.57 & $6.1\times10^{-7}$ & 4.47 & $1.3\times10^{-2}$ & 6.17 \\
        & 6511 & $1.1\times10^{21}$ & $1.2\times10^{-7}$ & 89.49& $8.5\times10^{-5}$ & 34.09& $6.1\times10^{-7}$ & 30.56& $4.9\times10^{-2}$ & 42.67 \\
\bottomrule
\end{tabular}
\end{table}


\subsection{Example 2: Helmholtz Equation in the Unit Sphere}
\begin{figure}
\begin{center}
\includegraphics[width=.32\textwidth,height=.35\textwidth]{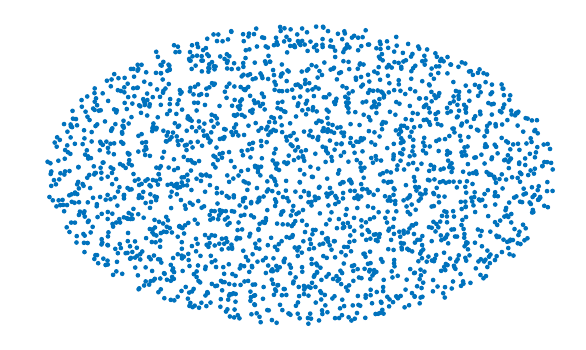}%
\includegraphics[width=.35\textwidth,height=.33\textwidth]{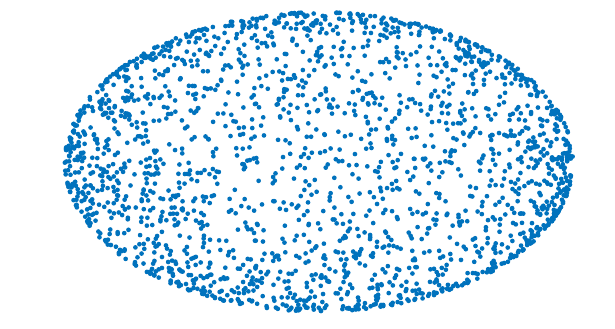}%
\includegraphics[width=.32\textwidth,height=.34\textwidth]{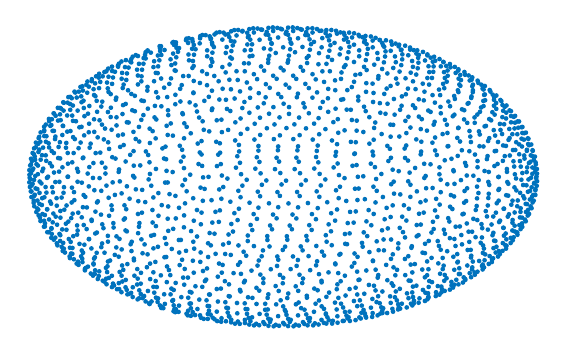}
\caption{Collocation points on the unit sphere for Example~2: (left) random distribution on the sphere, (middle) Halton sequence on the sphere, (right) uniform distribution.}
    \label{fig:points_sphere}
\end{center}
\end{figure}

\begin{figure}
\begin{center}
\includegraphics[width=.5\textwidth,height=.4\textwidth]{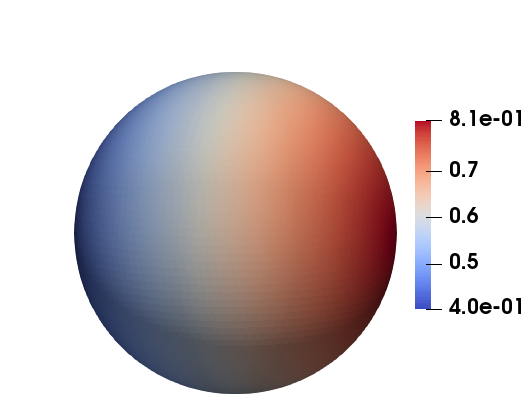}%
\includegraphics[width=.5\textwidth,height=.4\textwidth]{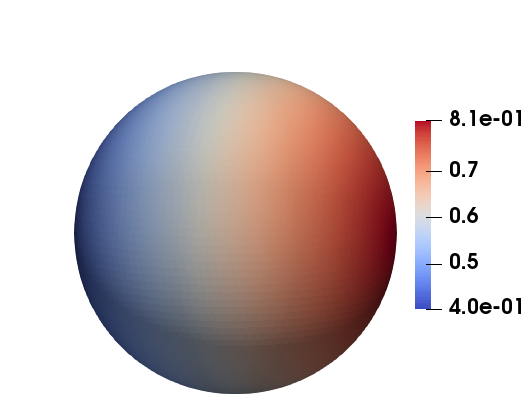}
\caption{Example~2 (Unit Sphere): (Left) geometry of the unit sphere domain with collocation points, and (Right) the approximate solution obtained by the HKT(180) method for $k=3$ with $N=2154$ uniformly distributed points (surface plot on the sphere). The solution is smooth and no spurious oscillations are visible, indicating effective regularization.}
    \label{fig:solution_sphere}
\end{center}
\end{figure}

In this example, we solve the Helmholtz problem \eqref{eq:Helmholtz} inside a unit sphere. The boundary of the domain is the unit sphere $x^2+y^2+z^2 = 1$, on which we impose Dirichlet boundary conditions derived from a known analytical solution. We choose the exact solution in the spherical domain to be of the same form as in Example~1:
\[
u(x,y,z) = \exp\!\Big(-\frac{(x-0.25)^2 + (y-0.25)^2 + z^2}{\sigma}\Big)\,,
\qquad \sigma = 20\,,
\] 
so that $u(x,y,z)$ is slightly off-centered within the sphere. Figure~\ref{fig:points_sphere} shows the collocation points on the sphere for a representative random, uniform, and Halton distribution. Figure~\ref{fig:solution_sphere} illustrates the geometry of the unit sphere and a surface plot of the approximate solution obtained by the HKT(180) method for $k=3$ with $N=2154$ points.

We again compare the performance of the Tikh-Reg, Ine-TSVD, HKT, and Reg-GMRES methods. Table~\ref{tabspheric} summarizes the results for Example~2 (unit sphere) at $k=3$. The trends are similar to those observed in Example~1. All methods produce accurate results, with the hybrid approach (HKT) yielding the smallest errors in most cases. Notably, the sphere geometry appears slightly better conditioned than the cube (as reflected by smaller 
$\kappa(A)$ values in some cases), leading to smaller overall errors. The Reg-GMRES method benefits from this and achieves very small errors (on the order of $10^{-6}$ or lower) for the larger point sets, though at the cost of longer runtimes.
Again, the hybrid HKT approach provides
the best compromise between accuracy and efficiency, fully consistent with
Theorem \ref{thm:hybrid}.
\begin{table}[ht]
\centering
\caption{Unit sphere, $k = 3$ (Example~2). Relative errors (Re) and CPU time (s). $\mathrm{Re}$ and $\kappa(A)$ are rounded.}
\label{tabspheric}
\tiny
\setlength{\tabcolsep}{5pt}
\renewcommand{\arraystretch}{1.15}
\rowcolors{4}{gray!4}{white}
\begin{tabular}{l c c cc cc cc cc}
\toprule
Point distribution & $N$ & $\kappa(A)$
& \multicolumn{2}{c}{Tikh-Rg} & \multicolumn{2}{c}{Ine-TSVD} & \multicolumn{2}{c}{HKT}
& \multicolumn{2}{c}{Reg-GMRES} \\
\cmidrule(lr){4-5}\cmidrule(lr){6-7}\cmidrule(lr){8-9}\cmidrule(lr){10-11}
 &  &  & Re & CPU & Re & CPU & Re & CPU & Re & CPU \\
\midrule
Random  & 359  & $8.9\times10^{17}$ & $6.4\times10^{-9}$ & 0.09 & $1.7\times10^{-7}$ & 0.09 & $3.5\times10^{-8}$ & 0.07 & $3.6\times10^{-1}$ & 0.29 \\
        & 2154 & $4.8\times10^{20}$ & $6.2\times10^{-11}$ & 7.23 & $1.5\times10^{-9}$ & 4.97 & $8.0\times10^{-11}$ & 4.81 & $1.1\times10^{0}$ & 4.64 \\
        & 6511 & $2.8\times10^{21}$ & $8.0\times10^{-11}$ & 92.47& $1.7\times10^{-9}$ & 32.76& $4.1\times10^{-10}$ & 33.56& $7.8\times10^{-1}$ & 51.19 \\
\addlinespace[2pt]
Uniform & 359  & $3.4\times10^{18}$ & $1.2\times10^{-4}$ & 0.18 & $3.5\times10^{-4}$ & 0.28 & $2.4\times10^{-4}$ & 0.10 & $5.8\times10^{-6}$ & 0.29 \\
        & 2154 & $3.6\times10^{19}$ & $1.2\times10^{-6}$ & 7.45 & $1.3\times10^{-6}$ & 5.60 & $1.5\times10^{-6}$ & 4.71 & $1.2\times10^{0}$ & 5.97 \\
        & 6511 & $2.4\times10^{20}$ & $3.3\times10^{-6}$ & 93.65& $3.1\times10^{-6}$ & 36.16& $3.1\times10^{-6}$ & 34.66& $2.2\times10^{0}$ & 44.88 \\
\addlinespace[2pt]
Halton  & 359  & $5.3\times10^{16}$ & $1.3\times10^{-4}$ & 0.19 & $6.0\times10^{-4}$ & 0.18 & $7.1\times10^{-4}$ & 0.11 & $1.1\times10^{0}$ & 0.31 \\
        & 2154 & $6.4\times10^{19}$ & $4.5\times10^{-5}$ & 7.46 & $7.4\times10^{-4}$ & 5.50 & $1.9\times10^{-4}$ & 4.81 & $9.7\times10^{-1}$ & 5.98 \\
        & 6511 & $2.7\times10^{20}$ & $6.4\times10^{-5}$ & 83.03& $2.4\times10^{-2}$ & 36.16& $8.6\times10^{-5}$ & 31.66& $1.4\times10^{0}$ & 56.67 \\
\bottomrule
\end{tabular}
\end{table}

\subsection{Example 3: Helmholtz Equation in a Complex Geometry}
\begin{figure}[htbp]
	\centering
	\includegraphics[width=0.50\textwidth]{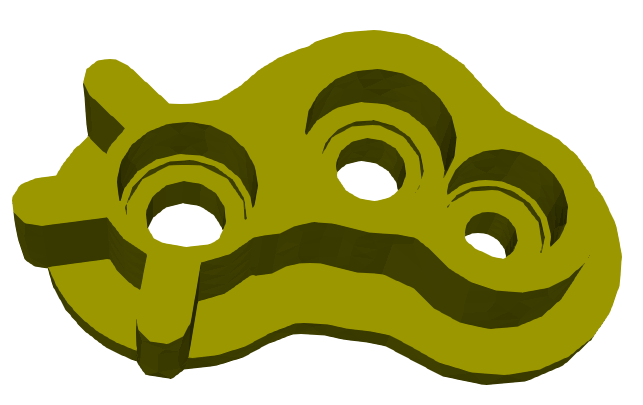}
	\caption{Geometry of the pump casing model (Example~3). This complex domain is used to test the proposed methods on an irregular three-dimensional geometry. Collocation points are distributed throughout the volume (surface not shown for clarity).}
    \label{fig:pump_geom}
\end{figure}

Finally, we test our methods on a complex three-dimensional geometry derived from a pump casing model(see Figure~\ref{fig:pump_geom}).  
The exact solution is again taken as a Gaussian field,
\[
u(x,y,z) = \exp\!\Big(-\frac{x^2 + y^2 + z^2}{\sigma}\Big), \qquad \sigma = 10,
\]  
with Dirichlet data prescribed on the entire boundary.

We use $N=1684$ and $N=6794$ collocation points distributed quasi-uniformly.  
For this large-scale case, direct TSVD or LU decomposition is impractical due to memory and cost.  
We therefore compare only HKT and Reg-GMRES.  

\begin{figure} [htbp]
	\begin{center}

		\subfloat[Exact solution ]{\includegraphics[height=1.5in]{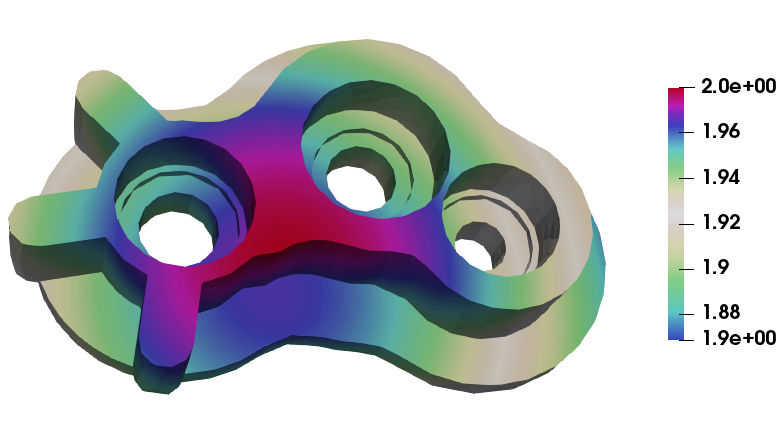}}
		\subfloat[Wave number $k=3$ ]{\includegraphics[height=1.5in]{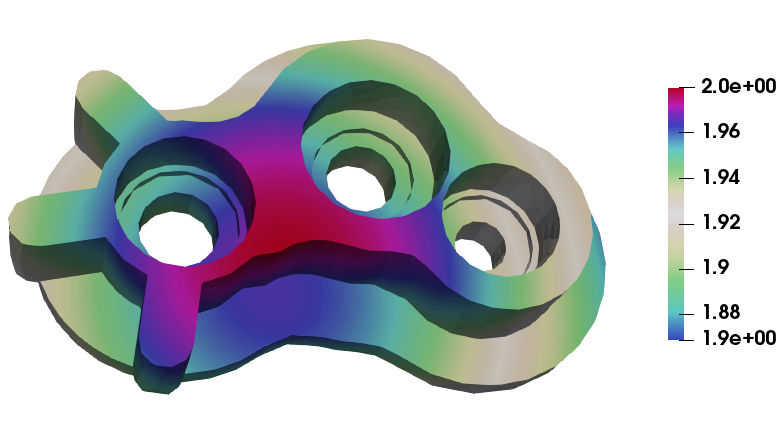}}\\
       \caption{Example~3 (Pump Casing): Exact solution (left) and approximate solution (right) for $k=3$ on the pump casing geometry, with $N=1684$ collocation points. The approximate solution is computed using the HKT method and is visually indistinguishable from the exact solution at the plotted scale.}
    \label{fig:pump_solution} 
	\end{center}
\end{figure}

\begin{figure}[htbp]
    \centering
    \includegraphics[width=0.6\textwidth,height=.6\textwidth]{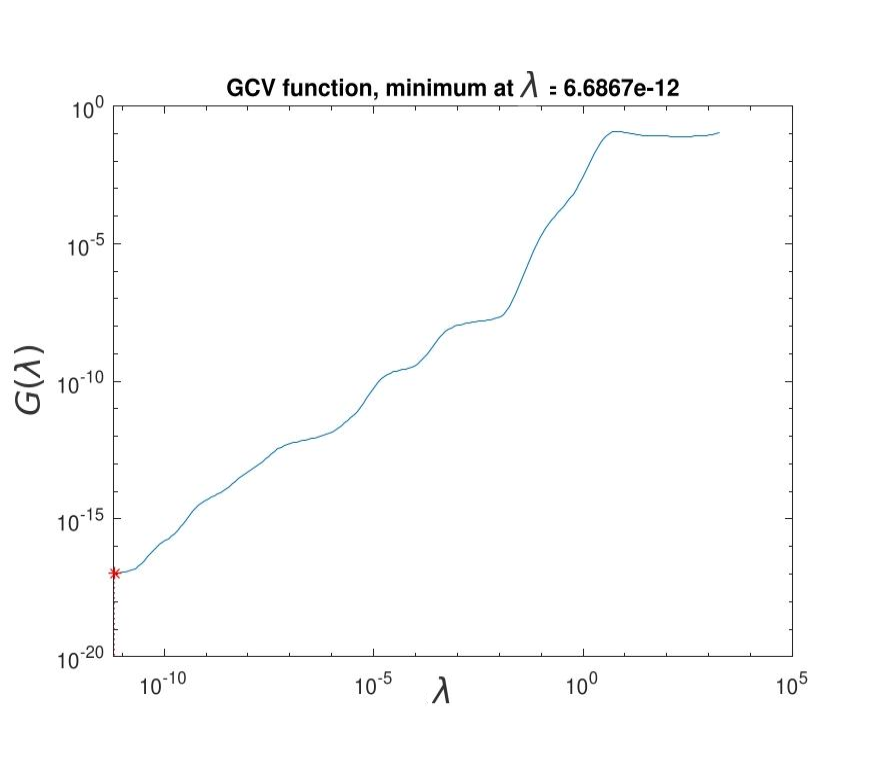}%
		\caption{GCV function for Example~3 ($k=3$, $N=1684$). The curve shows a well-defined minimum, indicated by the dot, at $\lambda \approx 6.7 \times 10^{-12}$. This optimal regularization parameter was used for the HKT solution.}
    \label{fig:gcv_curve}
\end{figure}

Figure~\ref{fig:pump_solution} illustrates the exact and HKT approximated solutions for $N=1684$, while Figure~\ref{fig:gcv_curve} plots the GCV function, showing a clear minimum at $\lambda \approx 6.7\times 10^{-12}$, which is used as the regularization parameter.  
Results are summarized in Table~\ref{tab:pump_results}.

HKT achieves accuracy comparable to a full TSVD (when feasible) but at drastically reduced cost.  
For $N=6794$, HKT attains an error of order $10^{-6}$ with manageable runtime, while Reg-GMRES converges much more slowly and delivers higher errors unless carefully stopped.  
This confirms the scalability of HKT to large, complex geometries.
\begin{table}[ht]
\centering
\caption{Pump-casing geometry, $k=3$ (Example~3). Relative error (Re) and CPU time (s). Re and $\kappa(A)$ are rounded}
\label{tab:pump_results}
\scriptsize
\setlength{\tabcolsep}{6pt}
\renewcommand{\arraystretch}{1.15}
\rowcolors{3}{gray!4}{white}
\begin{tabular}{c c cc cc}
\toprule
$N$ & $\kappa(A)$ & \multicolumn{2}{c}{HKT} & \multicolumn{2}{c}{Reg-GMRES} \\
\cmidrule(lr){3-4}\cmidrule(lr){5-6}
 &  & Re & CPU & Re & CPU \\
\midrule
1684 & $2\times 10^{20}$ & $2.05\times 10^{-6}$         & 10 & $2.99\times 10^{-1}$           & 12\\
6974 & $4\times 10^{20}$             & $8.37\times 10^{-6}$ & $650.20$              & $5.16\times 10^{-5}$ & $661.21$ \\
\bottomrule
\end{tabular}
\end{table}


\section{Conclusion}
\label{sec:conclusion}
This paper addressed a central obstacle to applying multiquadric (MQ) radial basis function collocation to the three-dimensional Helmholtz equation: the dense and ill-conditioned linear systems that arise from globally supported bases. We developed a solution framework that combines classical regularization with low-cost Krylov projections and demonstrated that it restores stability without sacrificing the approximation power that motivates meshless discretizations. Existence and uniqueness follow from standard arguments, while convergence and smoothing were established under a discrete Picard condition, clarifying how each method controls noise-amplifying modes.

Numerical experiments on a unit cube, a unit sphere, and a realistic industrial pump-casing geometry confirmed the effectiveness of this approach.  The results consistently showed that the hybrid scheme delivers the best accuracy–time balance, matching or surpassing the quality of full Tikhonov at substantially lower cost; the inexpensive truncated approach is the fastest when only the leading features of the solution are required; and the classical Tikhonov baseline remains competitive when the SVD information is readily available.

The proposed pipeline project, regularize, and lift offers a practical path to scalable, stable MQ-RBF methods for Helmholtz problems on complex three-dimensional domains.  Future work will target higher frequencies and larger node sets, improved automation of parameter and rank selection, and integration with localization strategies (e.g., domain decomposition or RBF-FD) and preconditioning. Extensions to variable-coefficient Helmholtz problems and coupled multi-physics settings are natural next steps, as the same projection-regularization principles apply.
\bibliographystyle{elsarticle-num}
\bibliography{ref}

\end{document}